# On Davenport and Heilbronn-Type of Functions

L. Ferry[1], D. Ghisa[2*] and F. A. Muscutar[3]

[1]*The Green Weasenham, Kings Lynn, Norfolk, PE32 2TD, United Kingdom*

[2]*York University, Glendon College, 2275 Bayview Avenue, Toronto, M4N 3M6, Canada*

[3]*Lorain CCC, 1005 Abbe Road, Eliria, OH44035, USA*

Opinion Article

## Abstract

A correction is brought to the opinion expressed in a previous note published in this journal that the off critical line points indicated as being non trivial zeros of Davenport and Heilbronn function are affected of approximation errors and illustrations are presented which enforce the conclusion that they are true zeros. It is shown also that linear combinations of L-functions satisfying the same Riemann-type of functional equation do not offer counterexamples to RH, contrary to a largely accepted position.

## 1. Introduction

In 1934 Potter and Titchmarsh [1], dealing with what is known today as Davenport and Heilbronn function, have shown that a certain linear combination of two Dirichlet L-functions corresponding to complex conjugate Dirichlet characters modulo 5 satisfies a Riemann-type of functional equation.

Since the corresponding Dirichlet series cannot be represented as an Euler product, they suspected that this function might have off critical line non trivial zeros.

They even thought they had identified two such zeros, yet they acknowledged that "the calculations are very cumbrous, and can hardly be considered conclusive".

However, after almost 60 years Spira's [2] calculations produced some more such zeros, which have found an undisputed place in the literature (see [3], [4], [5]), despite of a certain apparent ambiguity regarding the symmetric points with respect to the critical line. In [6] and [7] some more off critical line zeros of that function have been indicated and other functions obtained by a similar construction have been shown as possessing off critical line zeros.

By studying these functions, we found an incongruity, with which we will deal in section 3, and we had the impression that something similar must have taken place in Spira's calculations.

Our geometric approach contradicted the existence of such zeros, since it was showing no sign of symmetry with respect to the critical line.

The only explanation we had at hand was that some errors of approximation produced false off critical line zeros.

We made known our findings in [8]. However later we discovered that inadvertently one of the Dirichlet characters we were supposed to use in the construction of the Davenport and Heilbronn function was wrong, hence our function did not satisfy a Riemann-type of functional equation and, normally this fact produced the lack of symmetry.

The only way to interpret this lack of symmetry was to invoke the effect of approximation errors of the coefficients used in the respective linear combination.

After correcting that mistake and zooming on the region where Spira's zeros were located, the graphics have shown clearly not only those zeros, but also the symmetric ones with respect to the critical line.

We are now able to describe exactly what happens geometrically with the fundamental domains containing those zeros.

Although the embracing phenomenon concerns only curves $\Gamma_{k,j}$, $j \neq 0$, something similar happens here with $\Gamma_{k,0}$ and a component of the pre-image of a ray in a very close position to the negative real half axis, so that these two are embracing a curve $\Gamma_{k,1}$ or $\Gamma_{k,-1}$.

The effect on the fundamental domain associated with the embraced curve is that it becomes bounded to the right. It was impossible to imagine that such a configuration can exist before discovering it in the case of Davenport and Heilbronn function.

Fig. 1 below is illustrating three couples of Spira's zeros.

All the zeros at the right of the critical line are situated on curves $\Gamma_{k,0}$, as expected, while the zeros at the left of the critical line are situated one on a curve $\Gamma_{k,1}$ and two on curves $\Gamma_{k,-1}$.

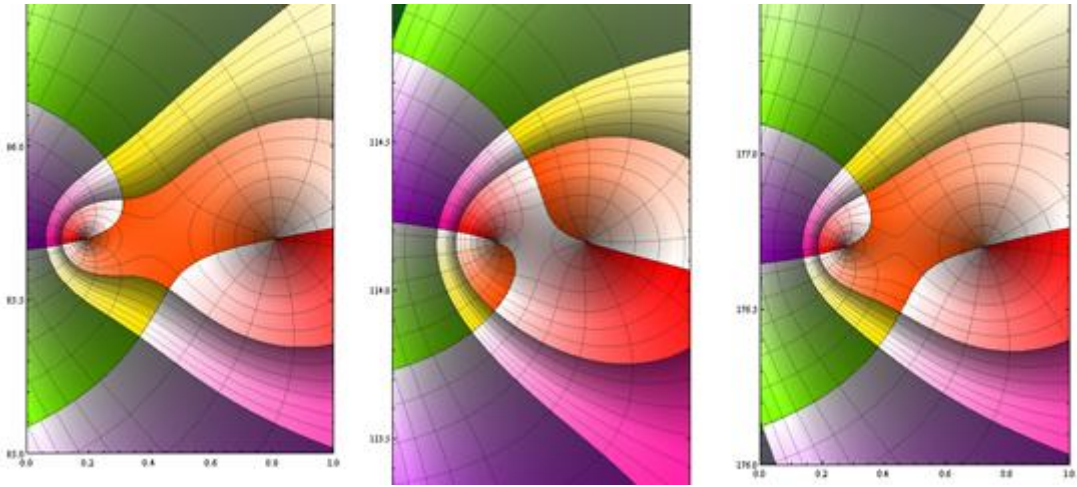

Fig.1 Zeros of Davenport and Heilbronn function
symmetric with respect to the critical line

In Fig. 2 we indicate the fundamental domain associated to an embraced curve $\Gamma_{k,-1}$. The slit corresponding to the respective fundamental domain is a ray $L$ starting at the image of the branch point (the zero of the derivative). The position of this zero is in accord with Speiser's theorem for Riemann Zeta function, hinting that the respective theorem might admit a generalization. It is approximately at 0.45+176.7i. The fact that it has the same imaginary part as the two zeros of the function is not a simple coincidence.

Indeed we have proved in [14] that if $f(s)$ satisfies a Riemann-type of functional equation, but it does not satisfy RH, then for every two distinct non trivial zeros $s_1 = \sigma + it$ and $s_2 = 1 - \sigma + it$ of $f(s)$ there is a zero $s_0$ of $f'(s)$ located on the interval determined by $s_1$ and $s_2$. It is the unique example encountered until now where the interval $(1, +\infty)$ of the real axis is not a part of the slit bounding the image of a fundamental domain.

Can such a configuration appear for a Dirichlet L-function? The answer is negative and a more general class of functions has been investigated from this point of view in [9], [10].

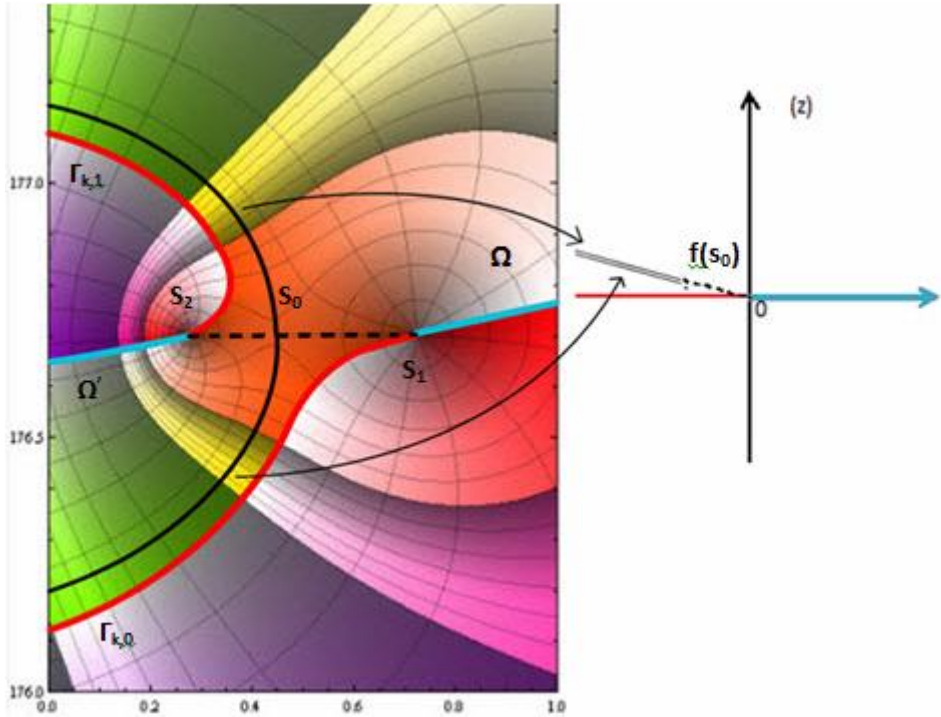

Fig. 2 Fundamental domains containing symmetric zeros
with respect to the critical line

### 2. Davenport and Heilbronn-Type of Function

In this section we take a new look at the Davenport and Heilbronn function. It has been obtained by analytic continuation to the whole complex plane of the Dirichlet series

$$f(s) = 1 + \frac{\tan\theta}{2^s} - \frac{\tan\theta}{3^s} - \frac{1}{4^s} + \frac{0}{5^s} + \frac{1}{6^s} + \ldots \qquad (1)$$

It has been proved that for $\tan^2\theta = [\sqrt{2} - \sqrt{1 + \frac{1}{\sqrt{5}}}] / [\sqrt{2} + \sqrt{1 + \frac{1}{\sqrt{5}}}]$, i.e. $\tan\theta = 0.284079\ldots$ this function can be written as:

$$f(s) = \tfrac{1}{2}\sec\theta[e^{-i\theta}L(s;\chi) + e^{i\theta}L(s;\overline{\chi})] \qquad (2)$$

where $\chi$ denotes the Dirichlet character modulo 5 for which $\chi(2) = i$.
This function satisfies the Riemann-type of functional equation

$$f(s) = 2^s\pi^{s-1}5^{(1/2)-s}\Gamma(1-s)\cos\tfrac{\pi s}{2}f(1-s) \qquad (3)$$

Since $f(s)$ is real for real $s$, if $f(\sigma_0 + it_0) = 0$ then necessarily $f(\sigma_0 - it_0) = 0$ and if $\cos[\pi(\sigma_0 + it_0)/2] \neq 0$, then due to (3), we have also $f(1 - \sigma_0 + it_0) = 0$, in other words the non trivial zeros of the function $f(s)$ are two by two symmetric with respect to the critical line.

In Fig. 1 the zeros from the right of the critical line are three of those indicated by Spira.

We are reluctant to call these zeros counterexamples to the Riemann Hypothesis (RH) since the function $f(s)$ is not a Dirichlet L-function for which the generalized RH has been stated, neither does it belong to the Selberg class for which the Grand RH is expected to be true.

They are simply an illustration of the fact that a Riemann-type of functional equation implies this symmetry of some non trivial zeros with respect to the critical line. Functions of this type can be constructed for an arbitrary modulus $q$ possessing complex conjugate characters $\chi$ and $\bar{\chi}$.

It is known (see [3], Corollary 10.9) that if $\chi$ is a primitive Dirichlet character modulo $q$, then the Dirichlet L-function

$$L(s;\chi) = \Sigma_{n=1}^{\infty} \chi(n)n^{-s} \qquad (4)$$

satisfies a Riemann-type of functional equation of the form:

$$L(s;\chi) = \epsilon(\chi)W(s)L(1-s;\bar{\chi}) \qquad (5)$$

where $W(s) = 2^s q^{(1/2)-s}\pi^{s-1}\Gamma(1-s)\sin\frac{\pi}{2}(s+\kappa)$, $\kappa = 0$ if $\chi(-1) = 1$ ($\chi$ is even), $\kappa = 1$ if $\chi(-1) = -1$ ($\chi$ is odd) and

$$\epsilon(\chi) = \tau(\chi) / i^\kappa \sqrt{q}, \quad \tau(\chi) = \Sigma_{k=1}^{q}\chi(k)Exp\{2k\pi i/q\} \qquad (6)$$

**Theorem 1**. *If for an arbitrary modulus $q$ there are complex conjugate primitive Dirichlet characters $\chi(\mod q)$ and $\bar{\chi}(\mod q)$, then a Davenport-Heilbronn type of function*

$$f(s) = \tfrac{1}{2}\{[L(s;\chi) + L(s;\bar{\chi})] + i\tan\theta[L(s;\chi) - L(s;\bar{\chi})]\} \qquad (7)$$

*can be built, which satisfies a Riemann-type of functional equation*

**Proof**: In principle all we need is to duplicate the computation from [11], which has been done for the particular value $q = 5$. The functions $L(s;\chi)$ and $L(s;\bar{\chi})$ verify functional equations of the form (5) for which $\kappa$ is the same since conjugate Dirichlet characters have the same parity. The corresponding $W(s)$ is the same too.

Obviously, $\tau(\chi) \neq \tau(\bar{\chi})$, since otherwise the row matrices corresponding to the two characters would be linearly dependent, which is excluded.

Finally, $\epsilon(\chi) \neq \epsilon(\bar{\chi}) = \overline{\epsilon(\chi)}$, i.e. $\epsilon(\chi) \neq \pm 1$. By an elementary computation and having in view (5) it can be shown that for $e^{-i\theta}\epsilon(\chi) = e^{i\theta}$ and thus $e^{i\theta}\overline{\epsilon(\chi)} = e^{-i\theta}$ we have:

$$f(s) = \tfrac{1}{2}\sec\theta[e^{-i\theta}L(s;\chi) + e^{i\theta}L(s;\bar{\chi})] =$$

$$\tfrac{1}{2}\sec\theta[e^{-i\theta}\epsilon(\chi)W(s)L(1-s;\bar{\chi}) + e^{i\theta}\overline{\epsilon(\chi)}W(s)L(1-s;\chi)] =$$

$$\frac{W(s)}{2}\sec\theta[e^{i\theta}L(1-s;\overline{\chi}) + e^{-i\theta}L(1-s;\chi)] = W(s)f(1-s) \quad (8).$$

**Example**: If $\chi = \chi_2(\mod 7)$, then $\theta$ can be determined such that

$f(s) = W(s)f(1-s)$, where $W(s) = 2^s 7^{(1/2)-s} \pi^{s-1} \Gamma(1-s)\cos\frac{\pi s}{2}$. Indeed, $\chi = \chi_2(\mod 7)$ is odd, hence $\kappa = 1$, thus $\sin\frac{\pi}{2}(s+\kappa) = \cos\frac{\pi s}{2}$ and $\tau(\chi) = \sum_{k=1}^{7}\chi(k)Exp\{2k\pi i/7\} =$

$e^{2\pi i/7} + \omega^2 e^{4\pi i/7} + \omega e^{6\pi i/7} - \omega e^{8\pi i/7} - \omega^2 e^{10\pi i/7} - e^{12\pi i/7}$, where $\omega = e^{\pi i/3}$. Having in view that $\omega^2 = -\overline{\omega}$ and that $e^{2(7-k)\pi i/7}$ and $e^{2k\pi i/7}$ are complex conjugate numbers, it can be easily found that $\tau(\chi) = 2i[\sin\frac{2\pi}{7} - \overline{\omega}\sin\frac{3\pi}{7} + \omega\sin\frac{\pi}{7}]$ and
$\epsilon(\chi) = \tau(\chi)/i\sqrt{7} = \frac{1}{\sqrt{7}}[(\sin\frac{\pi}{7} + 2\sin\frac{2\pi}{7} - \sin\frac{3\pi}{7}) + i\sqrt{3}(\sin\frac{\pi}{7} + \sin\frac{3\pi}{7})]$

The values of these trigonometric functions can be computed approximately and we find that $\tan\theta = [\text{Im}\,\epsilon(\chi)] / [\text{Re}\,\epsilon(\chi)] = 2.386161273...$

It is expected $f(s)$ to have off critical line non trivial zeros.

### 3. Linear Combinations of L-Functions Satisfying the Same Functional Equation

This topic made the object of the separate publication [12]. However, the remark which follows is more suited to an opinion article like this one. In [6], as well as in [4] the following example has been given, which seemed to close the discussion about the availability of counterexamples to GRH, since it has offered plenty of them.

Let $f_k(s)$, $k = 1,2$ be two L-functions satisfying the same Riemann-type of functional equation:

$$f_k(s) = W(s)\overline{f_k(1-\overline{s})} \quad (9)$$

Then, for an arbitrary complex number $s_0$, the function

$$f(s) = f_1(s_0)f_2(s) - f_2(s_0)f_1(s) \quad (10)$$

has obviously the zero $s_0$, which can be taken anywhere off the critical line. Then for the purpose of showing that a certain class of functions does not satisfy RH it looks like to be enough to find two linearly independent functions in the class satisfying the same Riemann-type of functional equation. It is implied in [6] and [4] that $f(s)$ satisfy also the functional equation (9).

We will show that this cannot be true. Indeed, by (9) and (10) we have respectively:

$$f(s) = f_1(s_0)W(s)\overline{f_2(1-\overline{s})} - f_2(s_0)W(s)\overline{f_1(1-\overline{s})} \quad (11)$$

$$W(s)\overline{f(1-\overline{s})} = W(s)\overline{[f_1(s_0)f_2(1-\overline{s}) - f_2(s_0)f_1(1-\overline{s})]} \quad (12)$$

Thus, if one of $f_k(s_0)$ is not real, the two expressions are different, hence $f(s)$ does not satisfy (9).

If both $f_1(s_0)$ and $f_2(s_0)$ are real and $\operatorname{Re} s_0 = 1/2$, there is nothing to prove. Suppose that $\operatorname{Re} s_0 \neq 1/2$, i.e. $s_0 \neq 1 - \overline{s_0}$. Since , by [13] the zeros of $f_k(s)$ are simple, if $f_1(s_0) = f_2(s_0) = 0$, hence $f_1(1 - \overline{s_0}) = f_2(1 - \overline{s_0}) = 0$, then the function $f_1(s) / f_2(s)$ has removable singularities at $s_0$ and $1 - \overline{s_0}$.

In this case and in the case where $f_2(s_0) \neq 0$, we have $f_1(s_0) / f_2(s_0) = f_1(1 - \overline{s_0}) / f_2(1 - \overline{s_0})$.

If $f_2(s_0) = 0$ and $f_1(s_0) \neq 0$, then $s_0$ and $1 - \overline{s_0}$ are simple poles of this function. Thus, $f_1(s) / f_2(s)$ takes the same value at $s_0$ and at $1 - \overline{s_0}$. Yet, for each one of $f_1(s)$ and $f_2(s)$ we can choose fundamental domains containing both $s_0$ and $1 - \overline{s_0}$. Indeed, if $f_k(s_0)$ and $f_k(1 - \overline{s_0})$, $k = 1, 2$ belong to $(-\infty, 1)$, any fundamental domain built as usual will contain these points. If some of these values are on the interval $(1, +\infty)$, it is always possible to slightly alter the boundaries of those domains such that the domains include those points.

Then the ratio $f_1(s) / f_2(s)$ which should be univalent in the intersection of those domains, takes the same value at $s_0$ and $1 - \overline{s_0}$, which is a contradiction.

The final conclusion is that this apparently simple counterexample to GRS is in fact not valid.

### 4. Conclusion

We bring a correction to a previous publication in this journal which assumed that approximation errors were responsible for some points appearing as off critical line zeros of Davenport and Heilbronn function.

Illustrations are presented showing that these points are indeed true zeros.

An argumentation is offered to the fact that L-functions satisfying the same Riemann-type of functional equation do not offer valid counterexamples to GRH, contrary to a largely admitted opinion.